# Calculus Without Limit Theory

*Jingzhong Zhang*
School of Computer Science and Educational Software
Guangzhou University
Guangzhou 510006, China
zjz2271@163.com

*Zengxiang Tong*
Department of Mathematical Sciences
Otterbein University
Westerville, OH 43081
ztong@otterbein.edu

**ABSTRCT**   This paper establishes calculus upon two physical facts: (1) an average velocity is always between two instantaneous velocities, and (2) the motion of an object is determined once its velocity has been determined. It directly defines derivative and definite integral on an ordered field, proves the fundamental theorem of calculus with no auxiliary conditions, easily reveals the common properties of derivatives, and obtains derivative formulas for elementary functions. Further discussion shows that for continuously differentiable functions, the new definitions are in accord with the traditional concepts.

## 1    View Derivative through an Inequality

Let $u < v$ be two real numbers. We have $2u < u + v < 2v$, in other words,

$$2u < \frac{v^2 - u^2}{v - u} < 2v, \qquad (1\text{-}1)$$

The above inequality seemingly indicates the function $2x$ is the derivative of the function $f(x) = x^2$.

Does it show a general pattern?

If it does, then the above may have revealed an important fact: one can evaluate derivatives without using limit process which had been tried by many mathematicians, including Newton and Lagrange, without success.

In China, the high school students are required of finding the derivatives of the functions $x^3$, $\sqrt{x}$, and $\frac{1}{x}$. Their derivatives can be similarly obtained as follows:

For $f(x) = x^3$, we have $\frac{f(v)-f(u)}{v-u} = \frac{v^2-u^2}{v-u} = u^2 + uv + v^2$, which yields, when $0 < u < v$, $3u^2 < \frac{f(v)-f(u)}{v-u} < 3v^2$, indicating that the derivative of $f(x) = x^3$ is $3x^2$.

For $f(x) = \sqrt{x}$, we have $\frac{f(v)-f(u)}{v-u} = \frac{\sqrt{v}-\sqrt{u}}{v-u} = \frac{1}{\sqrt{v}+\sqrt{u}}$, which yields, when $0 < u < v$, $\frac{1}{2\sqrt{v}} < \frac{f(v)-f(u)}{v-u} < \frac{1}{2\sqrt{u}}$, indicating that the derivative of $f(x) = \sqrt{x}$ is $\frac{1}{2\sqrt{x}}$.

For $f(x) = \frac{1}{x}$, we have $\frac{f(v)-f(u)}{v-u} = \frac{\frac{1}{v}-\frac{1}{u}}{v-u} = -\frac{1}{uv}$, which yields, when $0 < u < v$, $-\frac{1}{u} < \frac{f(v)-f(u)}{v-u} < -\frac{1}{v^2}$, indicating that the derivative of $f(x) = \frac{1}{x}$ is $-\frac{1}{x^2}$.

Newton was not perfectly satisfied with his calculus, wishing to eliminate the vague "infinitesimal" from his fluxional calculus. He showed some thoughts about "limit" in his *Principia*. Historians believe that he had spent years to think over how to eliminate the vague "infinitesimal" and to perfect his theory without success.

Lagrange endeavored to establish calculus without infinitesimal or limit. He tried to use infinite series to overcome the dilemma of infinitesimal and limit. However, the convergence of a series cannot be explored without the concept of limit.

It has been an unshakable belief and tradition that calculus cannot be established without limit theory.

However, tradition is changeable. Inspired by the above discovery, we are going to show how to rigorously establish calculus without limit theory.

## 2    From Instantaneous Velocity to Difference-Quotient Control Function

To reveal the general pattern contained in (1-1) — (1-4), we explore the instantaneous velocity problem, a classical physical model of differential calculus.

Let $S = S(t)$ denote the location of a moving object at time $t$ along a straight line, and $V = V(t)$ the instantaneous velocity at the time $t$. Then, the average velocity of the object over a time period $[u, v]$ must be between instantaneous velocities at two moments in the time period. That is, there are $p$ and $q$ in $[u, v]$ such that the following inequality holds:

$$V(p) \le \frac{S(v)-S(u)}{v-u} \le V(q).$$

The abstraction of the characteristic of instantaneous velocity leads to the concept of Difference-Quotient Control Function.

**Definition 1**  Assume functions $F$ and $f$ be defined on a set $S$. If for any two points $u < v$ in $S$, there are $p$ and $q$ in the set $[u, v] \cap S$ such that

$$f(p) \le \frac{F(v)-F(u)}{v-u} \le f(q), \tag{2-1}$$

then the function $f$ is called a **difference-quotient control function** of $F$ on the set $S$. The above inequality is called **difference-quotient control inequality**. In cases with no confusion, difference-quotient control functions are called **control functions**.

Obviously, the functions $F(x)$ and $F(x) + C$ have same control functions.

Assume $f$ is a control function of $F$ on a set $S$, and $Q$ an interval contained in $S$. The following are true:

- If $f(x) \equiv 0$ in $Q$, then $F(x)$ is a constant in $Q$.
- If $f(x) \equiv C$ in $Q$, then $F(x)$ is a linear function in $Q$.
- If $f(x) > 0$ in $Q$, then $F(x)$ is increasing in $Q$.
- If $f(x) < 0$ in $Q$, then $F(x)$ is decreasing in $Q$.
- If $f(x)$ is increasing in $Q$, then $F(x)$ is convex down in $Q$.
- If $f(x)$ is decreasing in $Q$, then $F(x)$ is convex up in $Q$.

The above properties show that a control function contains plenty of information of the function controlled by it.

It is known from (1-1) to (1-4) that the functions $2x$, $3x^2$, $\frac{1}{2\sqrt{x}}$, and $-\frac{1}{x^2}$ are control functions of the functions $x^2$, $x^3$, $\sqrt{x}$, and $\frac{1}{x}$. Similarly, we can see that, for positive integer $n$, $g(x) = nx^{n-1}$ is a control function of the function $f(x) = x^n$ on the intervals $(-\infty, 0]$ and $[0, \infty)$; the function $\cos x$ is a control function of $\sin x$ on the intervals $[k\pi, k\pi + \frac{\pi}{2}]$ and $[k\pi - \frac{\pi}{2}, k\pi]$, see [10] and [12] for the proofs, where the control functions are named "B" functions.

How to show that $\cos x$ is a control function of $\sin x$ on the intervals $(-\infty, \infty)$? We have the following theorem:

**Proposition 1**  Assume $g(x)$ is a control function of $f(x)$ on an interval $I$, as well as on an interval $J$. Then the following hold:

(1) If an interval K is contained in I or J, then $g(x)$ is also a control function of $f(x)$ on an interval K.
(2) If $K = I \cup J$ and $I \cap J \neq \varphi$ then $g(x)$ is also a control function of $f(x)$ on an interval K.

**Proof**  The proof of (1) is left for readers. We are now proving (2).

According to the definition, we need to show if $[u, v]$ is a subinterval of K, then there are points p and q in $[u, v]$ such that

$$f(p) \leq \frac{F(v) - F(u)}{v - u} \leq f(q).$$

If $[u, v]$ is a subinterval of I or J, then the result is true from the proof of (1). Otherwise, without loss of generality, there is a point $a$ with $a \in I \cap J$, $[u, a] \subset I$ and $[a, v] \subset J$. Noting that the difference-quotient $\frac{F(v)-F(u)}{v-u}$ is between $\frac{F(a)-F(u)}{a-u}$ and $\frac{F(v)-F(a)}{v-a}$, and there are $p_1$ and $q_1$ in $[u, a]$, $p_2$ and $q_2$ in $[a, v]$, such that the value of $\frac{F(a)-F(u)}{a-u}$ is between $g(p_1)$ and $g(q_1)$, and the value of $\frac{F(v)-F(a)}{v-a}$ is between $g(p_2)$ and $g(q_2)$, we can see that the value of $\frac{F(v)-F(u)}{v-u}$ is between $g(p)$ and $g(q)$, where $g(p) = Min\{g(p_1), g(q_1), g(p_2), g(q_2)\}$ and $g(q) = Max\{g(p_1), g(q_1), g(p_2), g(q_2)\}$. The proof is completed.

**Example 1**     Show that $g(x) = 3x^2 + 2ax + b$ is a control function of $f(x) = x^3 + ax^2 + bx + c$ on the interval $(-\infty, \infty)$.

**Proof**     The difference-quotient of $f(x)$ on the interval $[u, v]$ is

$$D = \frac{F(v) - F(u)}{v - u} = v^2 + uv + u^2 + a(u + v) + b.$$

Its differences from $g(u)$ and $g(v)$ are

$$D - g(u) = v^2 - u^2 + uv - u^2 + a(v - u) = (v - u)(v + 2u + a),$$

$$g(v) - D = v^2 - u^2 + v^2 - uv + a(v - u) = (v - u)(2v + u + a).$$

If $v > u \geq -\frac{a}{3}$, the above two expressions are positive, i.e., $g(u) < D < g(v)$, meaning that $g(x)$ is a control function of $f(x)$ on the interval $[-\frac{a}{3}, \infty)$. If $u < v \leq -\frac{a}{3}$, the above two expressions are negative, i.e., $g(v) < D < g(u)$, meaning that $g(x)$ is a control function of $f(x)$ on the interval $(-\infty, -\frac{a}{3}]$. Therefore, $g(x)$ is a control function of $f(x)$ on the interval $(-\infty, \infty)$. The proof is completed.

We have considered only one ordinary difference-quotient control inequality and apply it to a few common functions. But we have been able to use them to easily reveal important properties of many functions. These properties are usually taught in the late stage of a calculus course. The following are two more examples.

**Example 2**     Find an approximation for $\sqrt{10}$ and estimate the error.

Using the control function of $\sqrt{x}$, we have $\frac{1}{2\sqrt{v}} < \frac{\sqrt{v} - \sqrt{u}}{v - u} < \frac{1}{2\sqrt{u}}$, yielding

$\frac{1}{2\sqrt{10}} < \sqrt{10} - \sqrt{9} < \frac{1}{2\sqrt{9}}$, i.e., $\frac{1}{2\sqrt{10}} - \frac{1}{6} < \sqrt{10} - (3 + \frac{1}{6}) < 0$. We can choose $(3 + \frac{1}{6})$ as an approximation for $\sqrt{10}$, and its error is smaller than

$$\left| \frac{1}{2\sqrt{10}} - \frac{1}{6} \right| = \frac{\sqrt{10} - 3}{6\sqrt{10}} = \frac{1}{6(3\sqrt{10} + 10)} < \frac{1}{114}.$$

**Example 3**   If an increasing function $g(x)$ is a control function of $f(x)$ on an interval $I$, then $f(x)$ is a convex function on $I$.

**Proof**  Let $[u, v] \subset I$. Since $g(x)$ is increasing, we have

$$f\left(\frac{u+v}{2}\right) - f(u) < g\left(\frac{u+v}{2}\right)\left(\frac{v-u}{2}\right) < f(v) - f\left(\frac{u+v}{2}\right).$$

It implies $f\left(\frac{u+v}{2}\right) < \frac{f(v)+f(u)}{2}$, that is, $f(x)$ is a convex function on $I$.

Control functions are amazing. However, they are not necessary derivatives. For example, the Dirichlet function $g(x) = \begin{cases} 0 & x \in Q \\ 1 & x \in R - Q \end{cases}$ is a control function of the function $f(x) = 0.5x$, but not its derivative. Therefore, "derivatives" should have more features than "control functions" do.

## 3  Definition of Derivative without Limit

What conditions should be added to a control function to make it a derivative? We should follow the rule of simplicity and naturalism. The purpose of finding velocity is to understand the motion process. Actually, velocity and motion process are mutually determined. Similarly, the purpose of introducing derivative is to understand the original (or the controlled) function. Therefore, derivative and original functions are supposed to closely related and mutually determined. Based on this principle, we define derivatives as follows:

**Definition 2**  Suppose $f(x)$ is a control function of $F(x)$ on an interval $Q$. If every function $G(x)$ with $f(x)$ as its control function on $Q$ must take the form $G(x) = F(x) + C$, where $C$ is a constant, then $f(x)$ is said to be the **derivative** of $F(x)$ on an interval $Q$. In this case, the function $F(x)$ is said to be **differentiable** and an **original function** of $f(x)$ on the interval $Q$.

In other words, the derivative of $F(x)$ is the control function belonging to only $F(x)$.

Our definition of derivative is not equivalent to the traditional one and Lax's one [1]. However, for the $C^1$ functions, including all elementary functions, all the three definitions are equivalent, except at a few points. With our definition, if a function is differentiable on a closed interval $[a, b]$, then its derivative function is bounded on $[a, b]$. This property hold for Lax differentiable functions, but not for the traditionally differentiable functions.

What features, which are natural and easy to be operated mathematically, should be added to a control function to make it a derivative? We will address this problem in section 6.

## 4  Integral Systems and Definite Integrals

We are going to introduce an axiomatic definition of definite integral. Refereeing [6], [7], and [9], the readers can see how this idea gradually evolved in the past ten years. The essence of this idea is the geometric area properties depicted for curved trapezoids.

**Definition 3**  Let $Q$ be an interval, $f: Q \to R$, and $S: Q \times Q \to R$. If the following properties hold

(i)  Additivity: $\forall (u, v, w) \in Q^3$,    $S(u,v) + S(v,w) = S(u,w)$

(ii)  Intermediate Value Property: $\forall\, u < v$ in $Q$, $\exists\, p, q$ in $[u, v]$ such that
$$f(p)(v - u) \leq S(u,v) \leq f(q)(v - u),$$

then $S(u, v)$ is called an **integral system** of $f(x)$ on $Q$. If $f(x)$ has a unique integral system $S(u, v)$ on $Q$, then $f(x)$ is said to be **integrable** on $Q$, and the value of $S(u, v)$ is called the **definite integral** of $f(x)$ on $[u, v]$, denoted $S(u, v) = \int_u^v f(x)dx$.

The inequality in (ii) is equivalent to
$$f(p) \leq \frac{S(u,v)}{v-u} \leq f(q),$$
which shows close relationship between integrals and control functions.

Lax [1] writes the Riemann integral of $f(x)$ on $Q$ as $I\,(f, Q)$, emphasizing that an integration is an operation, with its input a function plus an interval and its output a number. To determine the value of $I\,(f, Q)$, one needs only two properties:

(1) The additivity of $I\,(f, Q)$ with respect to the interval Q: If $Q_1$ and $Q_2$ are two disjoint subintervals of Q, then $I(f, Q_1 + Q_2) = I(f, Q_1) + I(f, Q_2)$.

(2) The boundedness of $I(f, Q)$ with respect to $f$: If $\forall x \in Q \ \ m \leq f(x) \leq M$, then $m|Q| \leq I(f, Q) \leq M|Q|$.

Readers can see the obvious similarity between Lax's and our ideas, though they are still different. Lax first assumes the Riemann integrability, then explore its properties. We use these two properties, plus the uniqueness, as axioms to build the theory of integrals. Further explorations show that our definition of integrals is equivalent to that of Riemann integrals. The interested readers can refer Chapter 17 of [10] for the detailed proof.

The following proposition shows the relationship between an integral system and a control function. Its proof is left for readers.

**Proposition 2** Suppose $S(u, v)$ is an integral system of $f(x)$ on $Q$. Let $c \in Q$, and $F(x) = S(c, x)$, then $f(x)$ is a control function of $F(x)$ on $Q$. Conversely, if $f(x)$ is a control function of $F(x)$ on $Q$, and $S(u, v) = F(v) - F(u)$, then $S(u, v)$ is an integral system of $f(x)$ on $Q$.

## 5  The Fundamental Theorem of Calculus

**Proposition 3 (Newton-Leibniz Formula)** If $f(x)$ is the derivative of $F(x)$ on an interval $Q$, then $\forall u \in Q, \ \forall v \in Q$,

$$\int_u^v f(x)dx = F(v) - F(u)$$

Vice versa, if $f(x)$ has a unique integral system $S(u, v) = \int_u^v f(x)dx$ on $Q$, for any fixed $u$ in $Q$, the function $F(x) = S(u, x) = \int_u^x f(t)dt$ has $f(x)$ as its derivative on $Q$.

**Proof** Assume $F'(x) = f(x)$ on $Q$. Since $f(x)$ is the control function of $F(x)$, we have from Proposition 2 that $S(u, v) = F(v) - F(u)$ is an integral system of $f(x)$ on $Q$. According to the definition of definite integral, we need to show $S(u, v)$ is the unique integral system, i.e., if $R(u, v)$ is also an integral system of $f(x)$ on $Q$, then $R(u, v) = S(u, v)$. In fact, let a point $a \in Q$ and $G(x) = R(a, x)$, then $f(x)$ is also a control function of $G(x)$ on $Q$. Since $f(x)$ is the derivative of $F(x)$, that is, $f(x)$ is the control function of only functions $F(x) + C$, where $C$ is a constant, implying $G(x) = F(x) + C$. Therefore, $R(u, v) = G(v) - G(u) = F(v) - F(u) = S(u, v)$, i.e., $S(u, v)$ is the unique integral system of $f(x)$ on $Q$.

Conversely, if $f(x)$ has a unique integral system $S(u, v) = \int_u^v f(x)dx$ on $Q$, then for any fixed $u$ in $Q$, we know from Proposition 2 that the function $F(x) = S(u, x) = \int_u^x f(t)dt$ has $f(x)$ as its control function on $Q$. It follows from the uniqueness of integral system that $f(x)$ is the control function of only functions $F(x) + C$, where $C$ is a constant, that is, $F'(x) = f(x)$. The proof is completed.

The above fundamental theorem is a little different from its counterpart in traditional calculus books. Our theorem does not require any continuity condition for $f(x)$, and it is a somewhat "if and only if" biconditional statement.

The property "uniqueness" plays an important role in the concepts of derivative and definite integral. However, it is difficult to verify the uniqueness in a general case. In the following, we are going to present some conditions which are sufficient for uniqueness and easy for being verified.

**Proposition 4**  If $f(x)$ is a function of bounded variation on an interval $Q$, and both $S(u,v)$ and $R(u,v)$ are integral systems of $f(x)$ on $Q$, then $S(u,v) = R(u,v)$.

**Proof**  Let $u = x_0 < x_1 < \cdots < x_n = v$ be points in $Q$. Since both $S(u,v)$ and $R(u,v)$ are integral systems of $f(x)$ on $Q$, there are $p_k$ and $q_k$ in $[x_{k-1}, x_k]$ for $k = 1, 2, \ldots, n$ such that

$$|S(u,v) - R(u.v)| \leq \sum_{k=1}^{n}(x_k - x_{k-1})|f(q_k) - f(p_k)|. \qquad (5\text{-}1)$$

And since $f(x)$ is a function of bounded variation on an interval $Q$, there is a constant M, independent of the choices of $u$, $v$, $x_k$, $p_k$, and $q_k$, such that

$$\sum_{k=1}^{n}|f(q_k) - f(p_k)| < M \qquad (5\text{-}2)$$

We are going to prove the proposition by contradiction. Assume that there are $u$ and $v$ in $Q$ with $|S(u,v) - R(u,v)| = d > 0$. Choose $\{x_k \,|\, k = 1, 2, \ldots, n\}$ equally dividing $[u, v]$. We have

$$d \leq \sum_{k=1}^{n}(x_k - x_{k-1})|f(q_k) - f(p_k)| < \frac{(v-u)M}{n} \qquad (5\text{-}3)$$

Selecting $n > \frac{(v-u)M}{d}$ will lead a contradiction. The proof is completed.

Since a piecewise monotonic function is of bounded variation, the following are immediate corollaries.

**Corollary 1** If $f(x)$ is piecewise monotonic on an interval $Q$, and both $S(u,v)$ and $R(u,v)$ are integral systems of $f(x)$ on $Q$, then $S(u,v) = R(u,v)$.

**Corollary 2** If $f(x)$ is a control function of $F(x)$ on an interval $Q$, and it is piecewise monotonic, then $F'(x) = f(x)$.

Therefore, we can conclude $(\sqrt{x})' = \frac{1}{2\sqrt{x}}$, $(\sin x)' = \cos x$, and $(x^n)' = nx^{n-1}$ for any positive integer $n$.

The following example shows our concept is not exactly same to the traditional derivative.

**Example 4** Let $f(x) = sgn(x)$ and $S(u, v) = |v| - |u|$, then $S(u, v)$ is the unique integral system of $f(x)$ on any interval $Q$.

This example shows $|x|' = sgn(x)$ on $(-\infty, \infty)$, while it is not the case for the traditional definition of derivative.

# 6 L-Derivative of a Function

Since piecewise monotonic functions have unique integral systems, the piecewise monotonic control functions must be derivatives. However the class of piecewise monotonic functions is not closed under fundamental operations. Lipschitz functions are of bounded variations on any finite interval, they are closed under fundamental operations, and they are general enough to meet the needs of calculus and its applications.

**Definition 4** Let $f(x)$ be a function defined on a closed interval $[a, b]$. If there is a positive number $M$ such that for any points $u$ and $v$ in $[a, b]$, the following inequality holds

$$|f(u) - f(v)| \leq M|u - v|$$

Then $f(x)$ is called a **Lipschitz function** on $[a, b]$, or a bounded difference-quotient function on $[a, b]$. The number $M$ is called a **Lipschitz constant** of $f(x)$ on $[a, b]$.

Obviously, Lipschitz functions are of bounded variations on any finite interval, and they are closed under fundamental operations addition, subtraction, multiplication, division, and composition. The following is an immediate consequence of Proposition 4.

**Corollary 3** If a Lipschitz function $f(x)$ is a control function of both functions $F(x)$ and $G(x)$ on $[a, b]$, then $F(x) - G(x)$ is a constant on $[a, b]$, that is, $f(x)$ is the derivative of $F(x)$ on $[a, b]$. In this case, $F(x)$ is often said to be **L-differentiable**, and $f(x)$ is often called the **L-derivative** of $F(x)$ on $[a, b]$.

**Proposition 5** If $g(x)$ is the L-derivative of $f(x)$ on $[a, b]$, $M$ is a Lipschitz constant of $g(x)$ on $[a, b]$, then for any $[u, v] \subset [a, b]$ and $s \in [u, v]$, the following inequality holds

$$\left|\frac{f(v)-f(u)}{v-u} - g(s)\right| \leq M|v - u|, \tag{6-1}$$

And equivalently,

$$|f(v) - f(u) - g(s)(v - u)| \leq M(v - u)^2 \tag{6-2}$$

**Proof** From the difference-quotient control inequality (2-1), there are points $p$ and $q$ in $[u, v]$ such that

$$g(p) \leq \frac{f(v)-f(u)}{v-u} \leq g(q).$$

We have

$$-M|v - u| \leq g(p) - g(s) \leq \frac{f(v)-f(u)}{v-u} - g(s) \leq g(q) - g(s) \leq M|v - u|,$$

yielding (6-1) and (6-2). The proof is completed.

We call the inequalities (6-1) and (6-2) as Linqun inequalities, because Lin used them to define derivatives, see [3], [5], [8].

**Proposition 6** If $g(x)$ and $f(x)$ are functions defined on $[a, b]$ and for any two distinct points $u$ and $v$ in $[a, b]$,

$$\left|\frac{f(v)-f(u)}{v-u} - g(u)\right| \leq M|v - u|,$$

then, $g(x)$ is a Lipschitz function and a control function of $f(x)$ on $[a, b]$.

**Proof** We first show that $g(x)$ is a Lipschitz function. Let $h = v - u$, we have from the given Linqu inequality that

$$-Mh^2 \leq f(v) - f(u) - g(u)h \leq Mh^2$$

and

$$-Mh^2 \leq f(u) - f(v) + g(v)h \leq Mh^2$$

Thus,

$$-2Mh^2 \leq (g(v) - g(u))h \leq 2Mh^2$$

implying

$$|g(v) - g(u)| \leq 2M|v - u|,$$

which means $g(x)$ is a Lipschitz function on $[a, b]$.

We are now going to prove that $g(x)$ is a control function of $f(x)$ on $[a, b]$, that is, for any $[u, v] \subset [a, b]$, there are two points $p$ and $q$ in $[u, v]$ such that

$$g(p) \leq \frac{f(v)-f(u)}{v-u} \leq g(q).$$

We first show $\exists p \in [u, v]$ such that $g(p) \leq \frac{f(v)-f(u)}{v-u}$. It is obviously true if for all $x \in (u, v]$ the difference-quotient $\frac{f(x)-f(u)}{x-u}$ is a constant. If they are not always a constant, then there is $[r, s] \subset [u, v]$ with

$$\frac{f(v)-f(u)}{v-u} - \frac{f(s)-f(r)}{s-r} = d > 0.$$

Divide the interval into $n$ equal subintervals with $h = \frac{s-r}{n}$. There must be at least one subinterval, say $[p, p+h]$, such that $\frac{f(p+h)-f(p)}{h} \leq \frac{f(s)-f(r)}{s-r}$. Choose $n$ making $Mh = \frac{M(s-r)}{n} < d$, and we have

$$g(p) \leq \frac{f(p+h)-f(p)}{h} + Mh < \frac{f(s)-f(r)}{s-r} + d = \frac{f(v)-f(u)}{v-u}.$$

Similarly, it can be shown that $\exists q \in [u, v]$ such that $(q) \geq \frac{f(v)-f(u)}{v-u}$. The proof is completed.

Therefore, we have reached an important result: the function $g(x)$ is a Lipschitz control function of $f(x)$ on $[a, b]$ if and only if there is a number $M > 0$, for any two distinct points $u$ and $v$ in $[a, b]$,

$$\left|\frac{f(v)-f(u)}{v-u} - g(u)\right| \leq M|v-u|.$$

The uniqueness of L-derivative and its operation rules can be obtained from the above sufficient-necessary condition. For the detailed proofs, please refer [10] and [12].

# 7    Concluding Remarks

Due to the limitation of the space, we could not address all important topics in calculus. The interested readers are strongly encouraged to refer [10], [12], [7], [9], [3], [5], and [8].

The information we want to convey is that limit theory is not a prerequisite for learning calculus and that high school students with knowledge of function can understand majority of calculus topics.

Calculus is mathematics of motion, and functions are mathematical model of motions. Therefore, our approach of teaching calculus will help students to clearly see the relationship between calculus concepts and physical worlds, will help students to learn faster, easier, and more effectively, and will help students to quickly grasp the most important techniques of calculus, and apply calculus to solve a lot mathematics and real life problems.

Of course, we do not mean that students should never study elegant limit theory. In fact, having been familiar with certain calculus knowledge, students may learn limit theory much better and much deeper. We hope our alternate approach will stimulate mathematics professors' creative thinking, thinking over the possibility of rebuilding calculus structure to enhance our students' learning quality.